# Mannheim Partner $D$-Curves in Minkowski 3-space $E_1^3$

**Mustafa Kazaz[a], H. Hüseyin Uğurlu[b], Mehmet Önder[a], Tanju Kahraman[a]**

[a]*Celal Bayar University, Department of Mathematics, Faculty of Arts and Sciences, , Manisa, Turkey.*
E-mails: mustafa.kazaz@bayar.edu.tr, mehmet.onder@bayar.edu.tr, tanju.kahraman@bayar.edu.tr
[b]*Gazi University, Gazi Faculty of Education, Department of Secondary Education Science and Mathematics Teaching, Mathematics Teaching Program, Ankara, Turkey.* E-mail: hugurlu@gazi.edu.tr

**Abstract**

In this paper, we give the definition, different types and characterizations of Mannheim partner $D$-curves in Minkowski 3-space $E_1^3$. We find the relations between the geodesic curvatures, the normal curvatures and the geodesic torsions of these associated curves. Furthermore, we show that the definition and the characterizations of Mannheim partner $D$-curves include those of Mannheim partner curves in some special cases in Minkowski 3-space $E_1^3$.

**1. Introduction**

In the study of the fundamental theory and the characterizations of space curves, the related curves for which there exist corresponding relations between the curves are very interesting and an important problem. The most fascinating examples of such curves are associated curves, the curves for which at the corresponding points of them one of the Frenet vectors of a curve coincides with the one of the Frenet vectors of the other curve. The well known of the associated curves is Bertrand curve which is characterized as a kind of corresponding relation between the two curves. The relation is that the principal normal of a curve is the principal normal of another curve i.e, the Bertrand curve is a curve which shares the normal line with another curve. Over years many mathematicians have studied on Bertrand curves in different spaces and consider the properties of these curves[1,2,3,5,6]. By considering the frame of the ruled surface, Ravani and Ku extended the notion of Bertrand curve to the ruled surfaces and named as Bertrand offsets[16]. The corresponding characterizations of the Bertrand offsets of timelike ruled surface were given by Kurnaz[10].

Furthermore, Bertrand curves are not only the example of associated curves. Recently, a new definition of the associated curves was given by Liu and Wang[11,22]. They called these new curves as Mannheim partner curves: Let $x$ and $x_1$ be two curves in the three dimensional Euclidean $E^3$. If there exists a corresponding relationship between the space curves $x$ and $x_1$ such that, at the corresponding points of the curves, the principal normal lines of $x$ coincides with the binormal lines of $x_1$, then $x$ is called a Mannheim curve, and $x_1$ is called a Mannheim partner curve of $x$. The pair $\{x, x_1\}$ is said to be a Mannheim pair. They showed that the curve $x_1(s_1)$ is the Mannheim partner curve of the curve $x(s)$ if and only if the curvature $\kappa_1$ and the torsion $\tau_1$ of $x_1(s_1)$ satisfy the following equation

$$\dot{\tau} = \frac{d\tau}{ds_1} = \frac{\kappa_1}{\lambda}(1 + \lambda^2 \tau_1^2)$$

for some non-zero constant $\lambda$. They also study the Mannheim curves in Minkowski 3-space[11,21]. Some different characterizations of Mannheim partner curves are given by Orbay and others[15]. Similar to the Bertrand offsets, they have also defined and characterized the Mannheim offsets of ruled surfaces[14]. The corresponding characterizations of the Mannheim offsets of timelike and spacelike ruled surfaces in Minkowski 3-space $E_1^3$ have given by Kazaz, Ugurlu and Onder[7,8].

The differential geometry of the curves fully lying on a surface in Minkowski 3-space $E_1^3$ is given by Ugurlu, Kocayigit and Topal[9,18,19,20]. They have given the Darboux frame of the curves according to the Lorentzian characters of surfaces and the curves.

In this paper we consider the notion of the Mannheim partner curve for the curves lying on the surfaces. We call these new associated curves as Mannheim partner $D$-curves and by using the Darboux frame of the curves we give the definition, different types and the characterizations of these curves in Minkowski 3-space $E_1^3$.

## 2. Preliminaries

The Minkowski 3-space $E_1^3$ is the real vector space $IR^3$ provided with the standart flat metric given by

$$\langle,\rangle = -dx_1^2 + dx_2^2 + dx_3^2$$

where $(x_1, x_2, x_3)$ is a rectangular coordinate system of $E_1^3$. An arbitrary vector $\vec{v} = (v_1, v_2, v_3)$ in $E_1^3$ can have one of three Lorentzian causal characters; it can be spacelike if $\langle \vec{v}, \vec{v} \rangle > 0$ or $\vec{v} = 0$, timelike if $\langle \vec{v}, \vec{v} \rangle < 0$ and null (lightlike) if $\langle \vec{v}, \vec{v} \rangle = 0$ and $\vec{v} \neq 0$. Similarly, an arbitrary curve $\vec{\alpha} = \vec{\alpha}(s)$ can locally be spacelike, timelike or null (lightlike), if all of its velocity vectors $\alpha'(s)$ are respectively spacelike, timelike or null (lightlike). We say that a timelike vector is future pointing or past pointing if the first compound of the vector is positive or negative, respectively. For any vectors $\vec{x} = (x_1, x_2, x_3)$ and $\vec{y} = (y_1, y_2, y_3)$ in $E_1^3$, in the meaning of Lorentz vector product of $\vec{x}$ and $\vec{y}$ is defined by

$$\vec{x} \times \vec{y} = \begin{vmatrix} e_1 & -e_2 & -e_3 \\ x_1 & x_2 & x_3 \\ y_1 & y_2 & y_3 \end{vmatrix} = (x_2 y_3 - x_3 y_2, x_1 y_3 - x_3 y_1, x_2 y_1 - x_1 y_2)$$

where

$$\delta_{ij} = \begin{cases} 1 & i = j, \\ 0 & i \neq j, \end{cases} \quad e_i = (\delta_{i1}, \delta_{i2}, \delta_{i3}) \text{ and } e_1 \times e_2 = -e_3, \; e_2 \times e_3 = e_1, \; e_3 \times e_1 = -e_2.$$

Denote by $\{\vec{T}, \vec{N}, \vec{B}\}$ the moving Frenet frame along the curve $\alpha(s)$ in the Minkowski space $E_1^3$. For an arbitrary spacelike curve $\alpha(s)$ in the space $E_1^3$, the following Frenet formulae are given,

$$\begin{bmatrix} \vec{T}' \\ \vec{N}' \\ \vec{B}' \end{bmatrix} = \begin{bmatrix} 0 & k_1 & 0 \\ -\varepsilon k_1 & 0 & k_2 \\ 0 & k_2 & 0 \end{bmatrix} \begin{bmatrix} \vec{T} \\ \vec{N} \\ \vec{B} \end{bmatrix},$$

where $\langle \vec{T}, \vec{T} \rangle = 1$, $\langle \vec{N}, \vec{N} \rangle = \varepsilon = \pm 1$, $\langle \vec{B}, \vec{B} \rangle = -\varepsilon$, $\langle \vec{T}, \vec{N} \rangle = \langle \vec{T}, \vec{B} \rangle = \langle \vec{N}, \vec{B} \rangle = 0$ and $k_1$ and $k_2$ are curvature and torsion of the spacelike curve $\alpha(s)$ respectively. Here, $\varepsilon$ determines the kind of spacelike curve $\alpha(s)$. If $\varepsilon = 1$, then $\alpha(s)$ is a spacelike curve with spacelike first principal normal $\vec{N}$ and timelike binormal $\vec{B}$. If $\varepsilon = -1$, then $\alpha(s)$ is a spacelike curve with timelike principal normal $\vec{N}$ and spacelike binormal $\vec{B}$. Furthermore, for a timelike curve $\alpha(s)$ in the space $E_1^3$, the following Frenet formulae are given in as follows,

$$\begin{bmatrix} \vec{T}' \\ \vec{N}' \\ \vec{B}' \end{bmatrix} = \begin{bmatrix} 0 & k_1 & 0 \\ k_1 & 0 & k_2 \\ 0 & -k_2 & 0 \end{bmatrix} \begin{bmatrix} \vec{T} \\ \vec{N} \\ \vec{B} \end{bmatrix}.$$

where $\langle \vec{T}, \vec{T} \rangle = -1$, $\langle \vec{N}, \vec{N} \rangle = \langle \vec{B}, \vec{B} \rangle = 1$, $\langle \vec{T}, \vec{N} \rangle = \langle \vec{T}, \vec{B} \rangle = \langle \vec{N}, \vec{B} \rangle = 0$ and $k_1$ and $k_2$ are curvature and torsion of the timelike curve $\alpha(s)$ respectively [18,19].

**Definition 2.1.** *i) Hyperbolic angle:* Let $\vec{x}$ and $\vec{y}$ be future pointing (or past pointing) timelike vectors in $IR_1^3$. Then there is a unique real number $\theta \geq 0$ such that $<\vec{x}, \vec{y}> = -|\vec{x}||\vec{y}|\cosh\theta$. This number is called the *hyperbolic angle* between the vectors $\vec{x}$ and $\vec{y}$.

*ii) Central angle:* Let $\vec{x}$ and $\vec{y}$ be spacelike vectors in $IR_1^3$ that span a timelike vector subspace. Then there is a unique real number $\theta \geq 0$ such that $<\vec{x}, \vec{y}> = |\vec{x}||\vec{y}|\cosh\theta$. This number is called the *central angle* between the vectors $\vec{x}$ and $\vec{y}$.

*iii) Spacelike angle:* Let $\vec{x}$ and $\vec{y}$ be spacelike vectors in $IR_1^3$ that span a spacelike vector subspace. Then there is a unique real number $\theta \geq 0$ such that $<\vec{x}, \vec{y}> = |\vec{x}||\vec{y}|\cos\theta$. This number is called the *spacelike angle* between the vectors $\vec{x}$ and $\vec{y}$.

*iv) Lorentzian timelike angle:* Let $\vec{x}$ be a spacelike vector and $\vec{y}$ be a timelike vector in $IR_1^3$. Then there is a unique real number $\theta \geq 0$ such that $<\vec{x}, \vec{y}> = |\vec{x}||\vec{y}|\sinh\theta$. This number is called the *Lorentzian timelike angle* between the vectors $\vec{x}$ and $\vec{y}$ [13,7,8].

**Definition 2.2.** A surface in the Minkowski 3-space $IR_1^3$ is called a timelike surface if the induced metric on the surface is a Lorentz metric and is called a spacelike surface if the induced metric on the surface is a positive definite Riemannian metric, i.e., the normal vector on the spacelike (timelike) surface is a timelike (spacelike) vector, [13,7,8].

**Lemma 2.1.** *In the Minkowski 3-space $IR_1^3$, the following properties are satisfied:*

*(i) Two timelike vectors are never orthogonal.*
*(ii) Two null vectors are orthogonal if and only if they are linearly dependent.*
*(iii) A timelike vector is never orthogonal to a null (lightlike) vector* [13,7,8].

## 3. Darboux Frame of a Curve Lying on a Surface in Minkowski 3-space $E_1^3$

Let $S$ be an oriented surface in three-dimensional Minkowski space $E_1^3$ and let consider a non-null curve $x(s)$ lying on $S$ fully. Since the curve $x(s)$ is also in space, there exists Frenet frame $\{\vec{T}, \vec{N}, \vec{B}\}$ at each points of the curve where $\vec{T}$ is unit tangent vector, $\vec{N}$ is principal normal vector and $\vec{B}$ is binormal vector, respectively.

Since the curve $x(s)$ lies on the surface $S$ there exists another frame of the curve $x(s)$ which is called Darboux frame and denoted by $\{\vec{T}, \vec{g}, \vec{n}\}$. In this frame $\vec{T}$ is the unit tangent of the curve, $\vec{n}$ is the unit normal of the surface $S$ and $\vec{g}$ is a unit vector given by $\vec{g} = \vec{n} \times \vec{T}$. Since the unit tangent $\vec{T}$ is common in both Frenet frame and Darboux frame, the vectors $\vec{N}$, $\vec{B}$, $\vec{g}$ and $\vec{n}$ lie on the same plane. Then, if the surface $S$ is an oriented timelike surface, the relations between these frames can be given as follows

| If the curve $x(s)$ is timelike. | If the curve $x(s)$ is spacelike. |
|---|---|

$$\begin{bmatrix} \vec{T} \\ \vec{g} \\ \vec{n} \end{bmatrix} = \begin{bmatrix} 1 & 0 & 0 \\ 0 & \cos\varphi & \sin\varphi \\ 0 & -\sin\varphi & \cos\varphi \end{bmatrix} \begin{bmatrix} \vec{T} \\ \vec{N} \\ \vec{B} \end{bmatrix}, \qquad \begin{bmatrix} \vec{T} \\ \vec{g} \\ \vec{n} \end{bmatrix} = \begin{bmatrix} 1 & 0 & 0 \\ 0 & \cosh\varphi & \sinh\varphi \\ 0 & \sinh\varphi & \cosh\varphi \end{bmatrix} \begin{bmatrix} \vec{T} \\ \vec{N} \\ \vec{B} \end{bmatrix}.$$

If the surface $S$ is an oriented spacelike surface, then the curve $x(s)$ lying on $S$ is a spacelike curve. So, the relations between the frames can be given as follows

$$\begin{bmatrix} \vec{T} \\ \vec{g} \\ \vec{n} \end{bmatrix} = \begin{bmatrix} 1 & 0 & 0 \\ 0 & \cosh\varphi & \sinh\varphi \\ 0 & \sinh\varphi & \cosh\varphi \end{bmatrix} \begin{bmatrix} \vec{T} \\ \vec{N} \\ \vec{B} \end{bmatrix}.$$

In all cases, $\varphi$ is the angle between the vectors $\vec{g}$ and $\vec{N}$.

According to the Lorentzian causal characters of the surface $S$ and the curve $x(s)$ lying on $S$, the derivative formulae of the Darboux frame can be changed as follows:

**i)** If the surface $S$ is a timelike surface, then the curve $x(s)$ lying on $S$ can be a spacelike or a timelike curve. Thus, the derivative formulae of the Darboux frame of $x(s)$ is given by

$$\begin{bmatrix} \dot{\vec{T}} \\ \dot{\vec{g}} \\ \dot{\vec{n}} \end{bmatrix} = \begin{bmatrix} 0 & k_g & -\varepsilon k_n \\ k_g & 0 & \varepsilon \tau_g \\ k_n & \tau_g & 0 \end{bmatrix} \begin{bmatrix} \vec{T} \\ \vec{g} \\ \vec{n} \end{bmatrix}, \quad \langle \vec{T}, \vec{T} \rangle = \varepsilon = \pm 1, \ \langle \vec{g}, \vec{g} \rangle = -\varepsilon, \ \langle \vec{n}, \vec{n} \rangle = 1. \tag{$1_a$}$$

**ii)** If the surface $S$ is a spacelike surface, then the curve $x(s)$ lying on $S$ is a spacelike curve. Thus, the derivative formulae of the Darboux frame of $x(s)$ is given by

$$\begin{bmatrix} \dot{\vec{T}} \\ \dot{\vec{g}} \\ \dot{\vec{n}} \end{bmatrix} = \begin{bmatrix} 0 & k_g & k_n \\ -k_g & 0 & \tau_g \\ k_n & \tau_g & 0 \end{bmatrix} \begin{bmatrix} \vec{T} \\ \vec{g} \\ \vec{n} \end{bmatrix}, \quad \langle \vec{T}, \vec{T} \rangle = 1, \ \langle \vec{g}, \vec{g} \rangle = 1, \ \langle \vec{n}, \vec{n} \rangle = -1. \tag{$1_b$}$$

In these formulae $k_g$, $k_n$ and $\tau_g$ are called the geodesic curvature, the normal curvature and the geodesic torsion, respectively. Here and in the following, we use "dot" to denote the derivative with respect to the arc length parameter of a curve.

The relations between geodesic curvature, normal curvature, geodesic torsion and $\kappa$, $\tau$ are given as follows

$$k_g = \kappa \cos\varphi, \ k_n = \kappa \sin\varphi, \ \tau_g = \tau + \frac{d\varphi}{ds}, \text{ if both } S \text{ and } x(s) \text{ are timelike or spacelike}, \tag{$2_a$}$$

$$k_g = \kappa \cosh\varphi, \ k_n = \kappa \sinh\varphi, \ \tau_g = \tau + \frac{d\varphi}{ds}, \text{ if } S \text{ is timelike and } x(s) \text{ is spacelike}. \tag{$2_b$}$$

(See [9,18,19]). Furthermore, the geodesic curvature $k_g$ and geodesic torsion $\tau_g$ of the curve $x(s)$ can be calculated as follows

$$k_g = \left\langle \frac{dx}{ds}, \frac{d^2x}{ds^2} \times n \right\rangle, \quad \tau_g = \left\langle \frac{dx}{ds}, n \times \frac{dn}{ds} \right\rangle \tag{3}$$

In the differential geometry of surfaces, for a curve $x(s)$ lying on a surface $S$ the followings are well-known

    **i)** $x(s)$ is a geodesic curve $\Leftrightarrow k_g = 0$,

    **ii)** $x(s)$ is an asymptotic line $\Leftrightarrow k_n = 0$,

    **iii)** $x(s)$ is a principal line $\Leftrightarrow \tau_g = 0$ [12].

Through every point of the surface passes a geodesic in every direction. A geodesic is uniquely determined by an initial point and tangent at that point. All straight lines on a surface are geodesics. Along all curved geodesics the principal normal coincides with the surface normal. Along asymptotic lines osculating planes and tangent planes coincide, along geodesics they are normal. Through a point of a nondevelopable surface pass two asymptotic lines which can be real or imaginary [17].

## 4. Mannheim Partner $D$-Curves in Minkowski 3-space $E_1^3$

In this section, by considering the Darboux frame, we define Mannheim partner $D$-curves and give the characterizations of these curves in Minkowski 3-space $E_1^3$.

**Definition 4.1.** Let $S$ and $S_1$ be oriented surfaces in three-dimensional Minkowski space $E_1^3$ and let consider the arc-length parameter curves $x(s)$ and $x_1(s_1)$ lying fully on $S$ and $S_1$, respectively. Denote the Darboux frames of $x(s)$ and $x_1(s_1)$ by $\{T, g, n\}$ and $\{T_1, g_1, n_1\}$, respectively. If there exists a corresponding relationship between the curves $x$ and $x_1$ such that, at the corresponding points of the curves, the Darboux frame element $g$ of $x$ coincides with the Darboux frame element $n_1$ of $x_1$, then $x$ is called a Mannheim $D$-curve, and $x_1$ is a Mannheim partner $D$-curve of $x$. Then, the pair $\{x, x_1\}$ is said to be a Mannheim $D$-pair. If there exist such curves lying on the oriented surfaces $S$ and $S_1$, respectively, we call the pair $\{S, S_1\}$ as Mannheim pair surfaces.

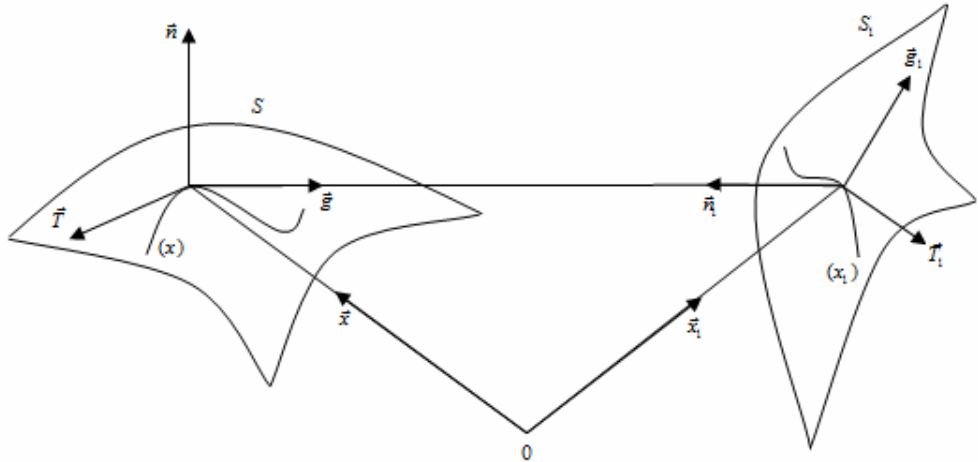

**Fig. 1** Mannheim partner $D$-curves

By considering the Lorentzian casual characters of the surfaces and the curves, from Definition 4.1, it is easily seen that there are five different types of the Mannheim $D$-curves in Minkowski 3-space. Let the pair $\{x, x_1\}$ be a Mannheim $D$-pair. Then according to the character of the surface $S$ we have the followings:

**Case 1.** The oriented surface $S$ is spacelike.
If both the surface $S$ and the curve $x(s)$ lying on $S$ are spacelike then, there are two cases; first one is that the surface $S_1$ is a timelike surface and the curve $x_1(s_1)$ fully lying on $S_1$ is spacelike. In this case we say that the pair $\{x, x_1\}$ is a Mannheim $D$-pair of the type 1. The second one is that both the surface $S_1$ and the curve $x_1(s_1)$ fully lying on $S_1$ are timelike. In this case we say that the pair $\{x, x_1\}$ is a Mannheim $D$-pair of the type 2.

**Case 2.** The oriented surface $S$ is timelike.
If the curve $x(s)$ lying on $S$ is a timelike curve then there are two cases; one is that both the surface $S_1$ and the curve $x_1(s_1)$ fully lying on $S_1$ are timelike. In this case we say that the pair $\{x, x_1\}$ is a Mannheim $D$-pair of the type 3. The other case is that the curve $x_1(s_1)$ fully lying on $S_1$ is a spacelike curve. Then the pair $\{x, x_1\}$ is a Mannheim $D$-pair of the type 4. If the curve $x(s)$ lying on $S$ is a spacelike curve then both the surface $S_1$ and the curve $x_1(s_1)$ fully lying on $S_1$ are spacelike. Then we say that the pair $\{x, x_1\}$ is a Mannheim $D$-pair of the type 5.

**Theorem 4.1.** *Let $S$ be an oriented surface and $x(s)$ be a Mannheim $D$-curve in $E_1^3$ with arc length parameter $s$ fully lying on $S$. If $S_1$ is another oriented surface and $x_1(s_1)$ is a curve with arc length parameter $s_1$ fully lying on $S_1$, then $x_1(s_1)$ is Mannheim partner $D$-curve of $x(s)$ if and only if the normal curvature $k_n$ of $x(s)$ and the geodesic curvature $k_{g_1}$, the normal curvature $k_{n_1}$ and the geodesic torsion $\tau_{g_1}$ of $x_1(s_1)$ satisfy the following equations*

*i) if the pair $\{x, x_1\}$ is a Mannheim $D$-pair of the type 1 or 3, then*

$$\dot{\tau}_{g_1} = \frac{1}{\lambda}\left[\left(\frac{(1+\lambda k_{n_1})^2 - \lambda^2 \tau_{g_1}^2}{(1+\lambda k_{n_1})}\right)\left(k_n \frac{1+\lambda k_{n_1}}{\cosh\theta} - k_{g_1}\right) + \frac{\lambda^2 \tau_{g_1} \dot{k}_{n_1}}{1+\lambda k_{n_1}}\right],$$

*ii) if the pair $\{x, x_1\}$ is a Mannheim $D$-pair of the type 2 or 4, then*

$$\dot{\tau}_{g_1} = \frac{1}{\lambda}\left[\left(\frac{(1+\lambda k_{n_1})^2 - \lambda^2 \tau_{g_1}^2}{(1+\lambda k_{n_1})}\right)\left(k_n \frac{1+\lambda k_{n_1}}{\sinh\theta} - k_{g_1}\right) + \frac{\lambda^2 \tau_{g_1} \dot{k}_{n_1}}{1+\lambda k_{n_1}}\right],$$

*iii) if the pair $\{x, x_1\}$ is a Mannheim $D$-pair of the type 5, we have*

$$\dot{\tau}_{g_1} = \frac{1}{\lambda}\left[\left(\frac{(1+\lambda k_{n_1})^2 + \lambda^2 \tau_{g_1}^2}{(1+\lambda k_{n_1})}\right)\left(-k_n \frac{1+\lambda k_{n_1}}{\cos\theta} - k_{g_1}\right) + \frac{\lambda^2 \tau_{g_1} \dot{k}_{n_1}}{1+\lambda k_{n_1}}\right],$$

*for some nonzero constants $\lambda$, where $\theta$ is the angle between the tangent vectors $T$ and $T_1$ at the corresponding points of the curves $x$ and $x_1$.*

**Proof: i)** Suppose that the pair $\{x, x_1\}$ is a Mannheim $D$-pair of the type 1. Denote the Darboux frames of $x(s)$ and $x_1(s_1)$ by $\{T, g, n\}$ and $\{T_1, g_1, n_1\}$, respectively. Then by the definition we can assume that

$$x(s_1) = x_1(s_1) + \lambda(s_1)n_1(s_1), \tag{4}$$

for some function $\lambda(s_1)$. By taking derivative of (4) with respect to $s_1$ and applying the Darboux formulas (1) we have

$$T\frac{ds}{ds_1} = (1+\lambda k_{n_1})T_1 + \dot{\lambda}n_1 + \lambda\tau_{g_1}g_1. \tag{5}$$

Since the direction of $n_1$ coincides with the direction of $g$, we get

$$\dot{\lambda}(s_1) = 0.$$

This means that $\lambda$ is a nonzero constant. Thus, the equality (5) can be written as follows

$$T \cdot \frac{ds}{ds_1} = (1 + \lambda k_{n_1})T_1 + \lambda \tau_{g_1} g_1. \tag{6}$$

On the other hand we have

$$T = \cosh\theta T_1 + \sinh\theta g_1, \tag{7}$$

where $\theta$ is the angle between the tangent vectors $T$ and $T_1$ at the corresponding points of $x$ and $x_1$. By differentiating this last equation with respect to $s_1$, we get

$$(k_g g + k_n n)\frac{ds}{ds_1} = (\dot\theta + k_{g_1})\sinh\theta T_1 + (\dot\theta + k_{g_1})\cosh\theta g_1 + (-k_{n_1}\cosh\theta + \tau_{g_1}\sinh\theta)n_1. \tag{8}$$

From this equation and the fact that

$$n = \sinh\theta T_1 + \cosh\theta g_1, \tag{9}$$

we get

$$(k_g g + k_n \sinh\theta T_1 + k_n \cosh\theta g_1)\frac{ds}{ds_1} = (\dot\theta + k_{g_1})\sinh\theta T_1 + (\dot\theta + k_{g_1})\cosh\theta g_1 \\ + (-k_{n_1}\cosh\theta + \tau_{g_1}\sinh\theta)n_1 \tag{10}$$

Since the direction of $n_1$ is coincident with $g$ we have

$$\dot\theta = k_n \frac{ds}{ds_1} - k_{g_1}. \tag{11}$$

From (6) and (7) and notice that $T_1$ is orthogonal to $g_1$ we obtain

$$\frac{ds}{ds_1} = \frac{1 + \lambda k_{n_1}}{\cosh\theta} = \frac{\lambda \tau_{g_1}}{\sinh\theta}. \tag{12}$$

Equality (12) gives us

$$\tanh\theta = \frac{\lambda \tau_{g_1}}{1 + \lambda k_{n_1}}. \tag{13}$$

By taking the derivative of this equation and applying (11) we get

$$\dot\tau_{g_1} = \frac{1}{\lambda}\left[\left(\frac{(1+\lambda k_{n_1})^2 - \lambda^2 \tau_{g_1}^2}{(1+\lambda k_{n_1})}\right)\left(k_n \frac{1+\lambda k_{n_1}}{\cosh\theta} - k_{g_1}\right) + \frac{\lambda^2 \tau_{g_1} \dot k_{n_1}}{1+\lambda k_{n_1}}\right], \tag{14}$$

that is desired.

Conversely, assume that the equation (14) holds for some nonzero constants $\lambda$. Then by using (12) and (13), (14) gives us

$$k_n \left(\frac{ds}{ds_1}\right)^3 = \lambda \dot\tau_{g_1}(1+\lambda k_{n_1}) - \lambda^2 \tau_{g_1} \dot k_{n_1} + \left((1+\lambda k_{n_1})^2 - \lambda^2 \tau_{g_1}^2\right)k_{g_1}. \tag{15}$$

Let define a curve

$$x(s_1) = x_1(s_1) + \lambda n_1(s_1). \tag{16}$$

where $\lambda$ is non-zero constant. We will prove that $x$ is a Mannheim $D$-curve and $x_1$ is the Mannheim partner $D$-curve of $x$. By taking the derivative of (16) with respect to $s_1$ twice, we get

$$\frac{ds}{ds_1}T = (1+\lambda k_{n_1})T_1 + \lambda \tau_{g_1} g_1, \tag{17}$$

and

$$(k_g g + k_n n)\left(\frac{ds}{ds_1}\right)^2 + T\frac{d^2s}{ds_1^2} = (\lambda \dot{k}_{n_1} + \lambda \tau_{g_1} k_{g_1})T_1 + \left((1+\lambda k_{n_1})k_{g_1} + \lambda \dot{\tau}_{g_1}\right)g_1 \qquad (18)$$
$$+ \left(-(1+\lambda k_{n_1})k_{n_1} + \lambda \tau_{g_1}^2\right)n_1$$

respectively. Taking the cross product of (17) with (18) we have

$$\left[k_g n + k_n g\right]\left(\frac{ds}{ds_1}\right)^3 = \left(\lambda \tau_{g_1} k_{n_1}(1+\lambda k_{n_1}) - \lambda^2 \tau_{g_1}^3\right)T_1 - \left[-(1+\lambda k_{n_1})^2 k_{n_1} + \lambda \tau_{g_1}^2 (1+\lambda k_{n_1})\right]g_1 \qquad (19)$$
$$+ \left[-k_{g_1}(1+\lambda k_{n_1})^2 - \lambda \dot{\tau}_{g_1}(1+\lambda k_{n_1}) + \lambda^2 \tau_{g_1} \dot{k}_{n_1} + \lambda^2 \tau_{g_1}^2 k_{g_1}\right]n_1$$

By substituting (15) in (19) we get

$$\left[k_g n + k_n g\right]\left(\frac{ds}{ds_1}\right)^3 = \left(\lambda \tau_{g_1} k_{n_1}(1+\lambda k_{n_1}) - \lambda^2 \tau_{g_1}^3\right)T_1$$
$$- \left(-k_{n_1}(1+\lambda k_{n_1})^2 + \lambda \tau_{g_1}^2 (1+\lambda k_{n_1})\right)g_1 \qquad (20)$$
$$- k_n \left(\frac{ds}{ds_1}\right)^3 n_1$$

Taking the cross product of (17) with (20) we have

$$\left[k_g g + k_n n\right]\left(\frac{ds}{ds_1}\right)^4 = k_n \left(\frac{ds}{ds_1}\right)^3 \lambda \tau_{g_1} T_1 + k_n \left(\frac{ds}{ds_1}\right)^3 (1+\lambda k_{n_1})g_1 \qquad (21)$$
$$+ \left((1+\lambda k_{n_1})^2 - \lambda^2 \tau_{g_1}^2\right)\left(\lambda \tau_{g_1}^2 - k_{n_1}(1+\lambda k_{n_1})\right)n_1$$

From (20) and (21) we have

$$(k_g^2 - k_n^2)\left(\frac{ds}{ds_1}\right)^4 \vec{n} = \left[k_g \frac{ds}{ds_1}\left(\lambda \tau_{g_1} k_{n_1}(1+\lambda k_{n_1}) - \lambda^2 \tau_{g_1}^3\right) - \lambda \tau_{g_1} k_n^2 \left(\frac{ds}{ds_1}\right)^3\right]\vec{T}_1$$
$$- \left[k_g \frac{ds}{ds_1}\left(-k_{n_1}(1+\lambda k_{n_1})^2 + \lambda \tau_{g_1}^2 (1+\lambda k_{n_1})\right) + (1+\lambda k_{n_1})k_n^2 \left(\frac{ds}{ds_1}\right)^3\right]\vec{g}_1 \quad (22)$$
$$- \left[k_n k_g \left(\frac{ds}{ds_1}\right)^4 + k_n \left((1+\lambda k_{n_1})^2 - \lambda^2 \tau_{g_1}^2\right)\left(\lambda \tau_{g_1}^2 - k_{n_1}(1+\lambda k_{n_1})\right)\right]\vec{n}_1$$

Furthermore, from (17) and (20) we get

$$\begin{cases} \left(\lambda^2 \tau_{g_1}^2 - (1+\lambda k_{n_1})^2\right) = \left(\frac{ds}{ds_1}\right)^2, \\ k_g \left(\frac{ds}{ds_1}\right)^2 - \lambda \tau_{g_1}^2 + k_{n_1}(1+\lambda k_{n_1}) = 0 \end{cases} \qquad (23)$$

respectively. Substituting (23) in (22) we obtain

$$(k_g^2 - k_n^2)\left(\frac{ds}{ds_1}\right)^4 \vec{n} = \left[k_g \frac{ds}{ds_1}\left(\lambda \tau_{g_1} k_{n_1}(1+\lambda k_{n_1}) - \lambda^2 \tau_{g_1}^3\right) - \lambda \tau_{g_1} k_n^2 \left(\frac{ds}{ds_1}\right)^3\right]\vec{T}_1$$
$$- \left[k_g \frac{ds}{ds_1}\left(-k_{n_1}(1+\lambda k_{n_1})^2 + \lambda \tau_{g_1}^2 (1+\lambda k_{n_1})\right) + (1+\lambda k_{n_1})k_n^2 \left(\frac{ds}{ds_1}\right)^3\right]\vec{g}_1 \qquad (24)$$

Equality (17) and (24) shows that the vectors $\vec{T}$ and $\vec{n}$ lie on the plane $sp\{\vec{T_1}, \vec{g_1}\}$. So, at the corresponding points of the curves, the Darboux frame element $\vec{g}$ of $x$ coincides with the Darboux frame element $\vec{n_1}$ of $x_1$, i.e, the curves $x$ and $x_1$ are Mannheim $D$-pair curves of the type 1. ∎

Let now give the characterizations of Mannheim partner $D$-curves in some special cases. Let the pair $\{x, x_1\}$ be a Mannheim $D$-pair of the type 1 or 3 in Minkowski 3-space $E_1^3$. Assume that $x(s)$ be an asymptotic Mannheim $D$-curve. Then, from (14) we have the following special cases:

**i)** Consider that $x_1(s_1)$ is a geodesic curve. Then $x_1(s_1)$ is Mannheim partner $D$-curve of $x(s)$ if and only if the following equation holds,
$$\dot{\tau}_{g_1} = -\frac{\lambda \tau_{g_1} \dot{k}_{n_1}}{1 - \lambda k_{n_1}}.$$

**ii)** Assume that $x_1(s_1)$ is also an asymptotic line. Then $x_1(s_1)$ is Mannheim partner $D$-curve of $x(s)$ if and only if the geodesic curvature $k_{g_1}$ and the geodesic torsion $\tau_{g_1}$ of $x_1(s_1)$ satisfy the following equation,
$$\lambda \dot{\tau}_{g_1} = (1 + \lambda^2 \tau_{g_1}^2) k_{g_1}.$$
In this case, the Frenet frame of the curve $x_1(s_1)$ coincides with its Darboux frame. From (2) we have $k_{g_1} = \kappa_1$ and $\tau_{g_1} = \tau_1$. So, in Minkowski 3-space the Mannheim partner $D$-curves become the Mannheim partner curves, i.e., if both $x(s)$ and $x_1(s_1)$ are asymptotic lines then, the definition and the characterizations of the Mannheim partner $D$-curves involve those of the Mannheim partner curves in Minkowski 3-space.

**iii)** If $x_1(s_1)$ is a principal line then $x_1(s_1)$ is Mannheim partner $D$-curve of $x(s)$ if and only if the geodesic curvature $k_{g_1} = 0$ i.e, $x_1(s_1)$ is also a geodesic curve or $k_{n_1} = -1/\lambda = const$.

The proofs of the statement (ii) and (iii) of Theorem 4.1 and the particular cases given above can be given by the same way of the proof of statement (i).

**Theorem 4.2.** *Let the pair $\{x, x_1\}$ be a Mannheim $D$-pair in Minkowski 3-space $E_1^3$. Then the relation between geodesic curvature $k_g$, geodesic torsion $\tau_g$ of $x(s)$ and the normal curvature $k_{n_1}$, the geodesic torsion $\tau_{g_1}$ of $x_1(s_1)$ is given as follows*

*i) if the pair $\{x, x_1\}$ is a Mannheim $D$-pair of the type 1, 3, 4 or 5 then*
$$k_g - k_{n_1} = \lambda(-k_g k_{n_1} + \tau_g \tau_{g_1}),$$
*ii) if the pair $\{x, x_1\}$ is a Mannheim $D$-pair of the type 2, then*
$$k_g + k_{n_1} = \lambda(-k_g k_{n_1} + \tau_g \tau_{g_1}).$$

**Proof: i)** Let $x(s)$ be a Mannheim $D$-curve and $x_1(s_1)$ be a Mannheim partner $D$-curve of $x(s)$ in Minkowski 3-space $E_1^3$ and the pair $\{x, x_1\}$ be of the type 1. Then by definition we can write
$$x_1(s_1) = x(s_1) - \lambda(s_1) g(s_1) \tag{25}$$
for some constants $\lambda$. By differentiating (25) with respect to $s_1$ we have

$$T_1 = (1 + \lambda k_g)\frac{ds}{ds_1}T - \lambda \tau_g \frac{ds}{ds_1}n \tag{26}$$

By the definition we have
$$T_1 = \cosh\theta T - \sinh\theta n. \tag{27}$$

From (26) and (27) we obtain
$$\cosh\theta = (1 + \lambda k_g)\frac{ds}{ds_1}, \quad \sinh\theta = \lambda \tau_g \frac{ds}{ds_1}. \tag{28}$$

Using (12) and (28) it is easily seen that
$$k_g - k_{n_1} = \lambda(-k_g k_{n_1} + \tau_g \tau_{g_1}). \qquad \blacksquare$$

From Theorem 2, we obtain the following special cases.

Let the pair $\{x, x_1\}$ be a Mannheim $D$-pair of the type 1 in Minkowski 3-space $E_1^3$. Then,

**i)** if $x_1$ is an asymptotic line, then
$$k_g = \lambda \tau_g \tau_{g_1}$$
**ii)** if $x_1$ is a principal line, then
$$k_g - k_{n_1} = -\lambda k_g k_{n_1}$$
**iii)** if $x$ is a geodesic curve, then
$$k_{n_1} = -\lambda \tau_g \tau_{g_1}$$
**iv)** if $x$ is a principal line then
$$k_g - k_{n_1} = -\lambda k_g k_{n_1}$$

The proof of the cases that the pair $\{x, x_1\}$ be a Mannheim $D$-pair of the type 2, 3, 4 or 5 can be given by a similar procedure used in the proof of the case that the pair $\{x, x_1\}$ is of the type 1.

**Theorem 4.3.** *Let $\{x, x_1\}$ be Mannheim $D$-pair of the type 1. Then the following relations hold:*

**i)** $k_{g_1} = k_n \dfrac{ds}{ds_1} - \dfrac{d\theta}{ds_1}$

**ii)** $\tau_g \dfrac{ds}{ds_1} = -k_{n_1} \sinh\theta + \tau_{g_1} \cosh\theta$

**iii)** $k_g \dfrac{ds}{ds_1} = -k_{n_1} \cosh\theta + \tau_{g_1} \sinh\theta$

**iv)** $\tau_{g_1} = \left(-k_g \sinh\theta + \tau_g \cosh\theta\right)\dfrac{ds}{ds_1}$

**Proof: i)** Since the pair $\{x, x_1\}$ is of the type 1, we have $\langle T, T_1 \rangle = \cosh\theta$. By differentiating this equality with respect to $s_1$ we have

$$\left\langle (k_g g + k_n n)\frac{ds}{ds_1}, T_1 \right\rangle + \left\langle T, k_{g_1} g_1 - k_{n_1} n_1 \right\rangle = \sinh\theta \frac{d\theta}{ds_1}.$$

Using the fact that the direction of $n_1$ coincides with the direction of $g$ and

$$T_1 = \cosh\theta T - \sinh\theta n,$$
$$g_1 = -\sinh\theta T + \cosh\theta n, \qquad (29)$$

we easily get that
$$k_{g_1} = k_n \frac{ds}{ds_1} - \frac{d\theta}{ds_1}.$$

**ii)** By definition we get $\langle n, n_1 \rangle = 0$. Differentiating this equality with respect to $s_1$ we have
$$\left\langle (k_n T + \tau_g g)\frac{ds}{ds_1}, n_1 \right\rangle + \left\langle n, k_{n_1} T_1 + \tau_{g_1} g_1 \right\rangle = 0.$$

By (29) we obtain
$$\tau_g \frac{ds}{ds_1} = -k_{n_1} \sinh\theta + \tau_{g_1} \cosh\theta.$$

**iii)** By differentiating the equation $\langle T, n_1 \rangle = 0$ with respect to $s_1$ we get
$$\left\langle (k_g g + k_n n)\frac{ds}{ds_1}, n_1 \right\rangle + \left\langle T, k_{n_1} T_1 + \tau_{g_1} g_1 \right\rangle = 0.$$

From (29) it follows that
$$k_g \frac{ds}{ds_1} = -k_{n_1} \cosh\theta + \tau_{g_1} \sinh\theta.$$

**iv)** By differentiating the equation $\langle g, g_1 \rangle = 0$ with respect to $s_1$ we obtain
$$\left\langle (-k_g T + \tau_g n)\frac{ds}{ds_1}, g_1 \right\rangle + \left\langle g, k_{g_1} T_1 + \tau_{g_1} n_1 \right\rangle = 0.$$

By considering (29) we get
$$\tau_{g_1} = \left(-k_g \sinh\theta + \tau_g \cosh\theta\right)\frac{ds}{ds_1}. \qquad \blacksquare$$

The statements of Theorem 4. 3 for the pairs $\{x, x_1\}$ of the type 2, 3, 4, and 5 can be given as follows and the proofs can be easily done by the same way of the case the pairs $\{x, x_1\}$ is of the type 1.

| For the pair $\{x, x_1\}$ of the type 2 | For the pair $\{x, x_1\}$ of the type 3 |
| --- | --- |
| **i)** $k_{g_1} = k_n \dfrac{ds}{ds_1} + \dfrac{d\theta}{ds_1}$ | **i)** $k_{g_1} = k_n \dfrac{ds}{ds_1} + \dfrac{d\theta}{ds_1}$ |
| **ii)** $\tau_g \dfrac{ds}{ds_1} = k_{n_1} \cosh\theta - \tau_{g_1} \sinh\theta$ | **ii)** $\tau_g \dfrac{ds}{ds_1} = k_{n_1} \sinh\theta - \tau_{g_1} \cosh\theta$ |
| **iii)** $k_g \dfrac{ds}{ds_1} = k_{n_1} \sinh\theta - \tau_{g_1} \cosh\theta$ | **iii)** $k_g \dfrac{ds}{ds_1} = k_{n_1} \cosh\theta - \tau_{g_1} \sinh\theta$ |
| **iv)** $\tau_{g_1} = \left(-k_g \cosh\theta + \tau_g \sinh\theta\right)\dfrac{ds}{ds_1}$ | **iv)** $\tau_{g_1} = \left(k_g \sinh\theta - \tau_g \cosh\theta\right)\dfrac{ds}{ds_1}$ |

| For the pair $\{x, x_1\}$ of the type 4 | For the pair $\{x, x_1\}$ of the type 5 |
|---|---|
| **i)** $k_{g_1} = k_n \dfrac{ds}{ds_1} + \dfrac{d\theta}{ds_1}$ | **i)** $k_{g_1} = -k_n \dfrac{ds}{ds_1} - \dfrac{d\theta}{ds_1}$ |
| **ii)** $\tau_g \dfrac{ds}{ds_1} = -k_{n_1} \cosh\theta + \tau_{g_1} \sinh\theta$ | **ii)** $\tau_g \dfrac{ds}{ds_1} = -k_{n_1} \cos\theta + \tau_{g_1} \sin\theta$ |
| **iii)** $k_g \dfrac{ds}{ds_1} = -k_{n_1} \sinh\theta + \tau_{g_1} \cosh\theta$ | **iii)** $k_g \dfrac{ds}{ds_1} = k_{n_1} \sin\theta + \tau_{g_1} \cos\theta$ |
| **iv)** $\tau_{g_1} = \left(k_g \cosh\theta - \tau_g \sinh\theta\right) \dfrac{ds}{ds_1}$ | **iv)** $\tau_{g_1} = \left(k_g \cos\theta + \tau_g \sin\theta\right) \dfrac{ds}{ds_1}$ |

Let now $x$ be a Mannheim $D$-curve and $x_1$ be a Mannheim partner $D$-curve of $x$ and the pair $\{x, x_1\}$ be of the type 1. From the first equation of (3) and by using the fact that $n_1$ is coincident with $g$ we have

$$k_{g_1} = \langle \dot{x}_1, \ddot{x}_1 \times n_1 \rangle = \langle \dot{x}_1, \ddot{x}_1 \times g \rangle$$

$$= k_n \left[(1+\lambda k_g)^2 - \lambda^2 \tau_g^2\right] \left(\dfrac{ds}{ds_1}\right)^3 + \left[-\lambda \dot{\tau}_g (1+\lambda k_g) + \lambda^2 \tau_g \dot{k}_g\right] \left(\dfrac{ds}{ds_1}\right)^2 \quad (30)$$

Then the relations between the geodesic curvature $k_{g_1}$ of $x_1(s_1)$ and the geodesic curvature $k_g$, the normal curvature $k_n$ and the geodesic torsion $\tau_g$ of $x(s)$ are given as follows,

If $k_g = 0$ then from (30) the geodesic curvature $k_{g_1}$ of $x_1(s_1)$ is

$$k_{g_1} = -\left(\dfrac{ds}{ds_1}\right)^3 (1+\lambda^2 \tau_g^2) k_n + \left(\dfrac{ds}{ds_1}\right)^2 \lambda \dot{\tau}_g. \quad (31)$$

If $k_n = 0$ then the relation between $k_g$, $\tau_g$ and $k_{g_1}$ is

$$k_{g_1} = \lambda \left(\dfrac{ds}{ds_1}\right)^2 \left(\dot{\tau}_g (1+\lambda k_g) - \lambda \tau_g \dot{k}_g\right). \quad (32)$$

If $\tau_g = 0$ then, for the geodesic curvature $k_{g_1}$, we have

$$k_{g_1} = -\left(\dfrac{ds}{ds_1}\right)^3 (1+\lambda k_g)^2 k_n. \quad (33)$$

From (31), (32) and (33) we give the following corollary.

***Corollary 4.1.*** *Let $x$ be a Mannheim $D$-curve and $x_1$ be a Mannheim partner $D$-curve of $x$ and the pair the pair $\{x, x_1\}$ be of the type 1. Then the relations between the geodesic curvature $k_{g_1}$ of $x_1(s_1)$ and the geodesic curvature, $k_g$, the normal curvature $k_n$ and the geodesic torsion $\tau_g$ of $x(s)$ are given as follows*

*i) If $x$ is a geodesic curve, then the geodesic curvature $k_{g_1}$ of $x_1(s_1)$ is*

$$k_{g_1} = k_n \left(1 - \lambda^2 \tau_g^2\right) \left(\dfrac{ds}{ds_1}\right)^3 - \lambda \dot{\tau}_g \left(\dfrac{ds}{ds_1}\right)^2$$

***ii)*** *If $x$ is an asymptotic line, then the equation of $k_{g_1}$ is*

$$k_{g_1} = \left[ -\lambda \dot{\tau}_g (1 + \lambda k_g) + \lambda^2 \tau_g \dot{k}_g \right] \left( \frac{ds}{ds_1} \right)^2$$

***iii)*** *If $x$ is a principal line, then the geodesic curvature $k_{g_1}$ of $x_1(s_1)$ is*

$$k_{g_1} = k_n (1 + \lambda k_g)^2 \left( \frac{ds}{ds_1} \right)^3.$$

If the pair $\{x, x_1\}$ is of the type 2, 3, 4 or 5 then the geodesic curvature of the curve $x_1(s_1)$ is given as follows,

| If the pair $\{x, x_1\}$ is of the type 2 | If the pair $\{x, x_1\}$ is of the type 3 |
|---|---|
| $k_{g_1} = k_n \left[ (1 + \lambda k_g)^2 + \lambda^2 \tau_g^2 \right] \left( \frac{ds}{ds_1} \right)^3$ $+ \left[ -\lambda \dot{\tau}_g (1 + \lambda k_g) + \lambda^2 \tau_g \dot{k}_g \right] \left( \frac{ds}{ds_1} \right)^2$ | $k_{g_1} = -k_n \left[ (1 - \lambda k_g)^2 + \lambda^2 \tau_g^2 \right] \left( \frac{ds}{ds_1} \right)^3$ $+ \left[ -\lambda \dot{\tau}_g (1 - \lambda k_g) - \lambda^2 \tau_g \dot{k}_g \right] \left( \frac{ds}{ds_1} \right)^2$ |

| If the pair $\{x, x_1\}$ is of the type 4 | If the pair $\{x, x_1\}$ is of the type 5 |
|---|---|
| $k_{g_1} = k_n \left[ (1 - \lambda k_g)^2 - \lambda^2 \tau_g^2 \right] \left( \frac{ds}{ds_1} \right)^3$ $+ \left[ \lambda \dot{\tau}_g (1 - \lambda k_g) + \lambda^2 \tau_g \dot{k}_g \right] \left( \frac{ds}{ds_1} \right)^2$ | $k_{g_1} = -k_n \left[ (1 - \lambda k_g)^2 + \lambda^2 \tau_g^2 \right] \left( \frac{ds}{ds_1} \right)^3$ $+ \left[ -\lambda \dot{\tau}_g (1 - \lambda k_g) - \lambda^2 \tau_g \dot{k}_g \right] \left( \frac{ds}{ds_1} \right)^2$ |

and the statements in Corollary 4.1 are obtained by the same way.

Similarly, if the pair $\{x, x_1\}$ is of the type 1, from the second equation of (3) and by using the fact that $g$ is coincident with $n_1$, the relation between the geodesic torsion $\tau_{g_1}$ of $x_1(s_1)$ and the geodesic torsion $\tau_g$ of $x(s)$ is given by

$$\tau_{g_1} = -\tau_g \left( \frac{ds}{ds_1} \right)^2. \tag{34}$$

Furthermore, by using (12), from (34) we have

$$\tau_g \tau_{g_1} = -\frac{\sinh^2 \theta}{\lambda^2}. \tag{35}$$

Then, from (34) and (35) we can give the following corollary.

***Corollary 4.2.*** *Let $x$ be a Mannheim $D$-curve and $x_1$ be a Mannheim partner $D$-curve of $x$ and $\{x, x_1\}$ be of the type 1. Then the relation between the geodesic torsion $\tau_{g_1}$ of $x_1(s_1)$ and the geodesic torsion $\tau_g$ of $x(s)$ is given by one of the followings,*

***i)*** $\tau_{g_1} = -\tau_g \left( \frac{ds}{ds_1} \right)^2$    *or*    ***ii)*** $\tau_g \tau_{g_1} = -\frac{\sinh^2 \theta}{\lambda^2}$

*and so, the Mannheim partner $D$-curve $x_1$ is a principal line when the Mannheim $D$-curve $x$ is a principal line.*

Similarly, from (12) and (34) we get

**iii)** $\dfrac{\tau_g}{\tau_{g_1}} = -\dfrac{\cosh^2 \theta}{(1+\lambda k_{n_1})^2}$

Then, if $x_1(s_1)$ is an asymptotic curve, i.e., $k_{n_1} = 0$, we have

$$\tau_g = -\cosh^2 \theta \tau_{g_1}. \tag{36}$$

From (36) we have the following corollary.

***Corollary 4.3.*** *Let $x$ be a Mannheim $D$-curve and $x_1$ be a Mannheim partner $D$-curve of $x$ and $\{x, x_1\}$ be of the type 1. If $x_1(s_1)$ is an asymptotic curve then the relation between the geodesic torsion $\tau_g$ of $x(s)$ and the geodesic torsion $\tau_{g_1}$ of $x_1(s_1)$ is given as follows,*

**iv)** $\tau_g = -\cosh^2 \theta \tau_{g_1}$

*where $\theta$ is the angle between the tangent vectors $T$ and $T_1$ at the corresponding points of $x$ and $x_1$.*

When the pair $\{x, x_1\}$ is of the type 2, 3, 4 or 5, then the relations which give the geodesic torsion $\tau_{g_1}$ of $x_1(s_1)$ are given as follows.

| For the pair $\{x, x_1\}$ of the type 2 | For the pair $\{x, x_1\}$ of the type 3 |
|---|---|
| **i)** $\tau_{g_1} = -\tau_g \left(\dfrac{ds}{ds_1}\right)^2$ | **i)** $\tau_{g_1} = \tau_g \left(\dfrac{ds}{ds_1}\right)^2$ |
| **ii)** $\tau_g \tau_{g_1} = -\dfrac{\cosh^2 \theta}{\lambda^2}$ | **ii)** $\tau_g \tau_{g_1} = \dfrac{\sinh^2 \theta}{\lambda^2}$ |
| **iii)** $\dfrac{\tau_g}{\tau_{g_1}} = -\dfrac{\sinh^2 \theta}{(1+\lambda k_{n_1})^2}$ | **iii)** $\dfrac{\tau_g}{\tau_{g_1}} = \dfrac{\cosh^2 \theta}{(1+\lambda k_{n_1})^2}$ |
| **iv)** $\tau_g = -\sinh^2 \theta \tau_{g_1}$, if $x_1(s_1)$ is an asymptotic curve. | **iv)** $\tau_g = \cosh^2 \theta \tau_{g_1}$, if $x_1(s_1)$ is an asymptotic curve. |
| For the pair $\{x, x_1\}$ of the type 4 | For the pair $\{x, x_1\}$ of the type 5 |
| **i)** $\tau_{g_1} = \tau_g \left(\dfrac{ds}{ds_1}\right)^2$ | **i)** $\tau_{g_1} = -\tau_g \left(\dfrac{ds}{ds_1}\right)^2$ |
| **ii)** $\tau_g \tau_{g_1} = \dfrac{\cosh^2 \theta}{\lambda^2}$ | **ii)** $\tau_g \tau_{g_1} = -\dfrac{\sin^2 \theta}{\lambda^2}$ |
| **iii)** $\dfrac{\tau_g}{\tau_{g_1}} = \dfrac{\sinh^2 \theta}{(1+\lambda k_{n_1})^2}$ | **iii)** $\dfrac{\tau_g}{\tau_{g_1}} = -\dfrac{\cos^2 \theta}{(1+\lambda k_{n_1})^2}$ |
| **iv)** $\tau_g = \sinh^2 \theta \tau_{g_1}$, if $x_1(s_1)$ is an asymptotic curve. | **iv)** $\tau_g = -\cos^2 \theta \tau_{g_1}$, if $x_1(s_1)$ is an asymptotic curve. |

## 4. Conclusions

In this paper, in Minkowski 3-space $E_1^3$, the definition and characterizations of Mannheim partner $D$-curves are given which is a new study of associated curves lying on surfaces. It is shown that in Minkowski 3-space $E_1^3$, the definition and the characterizations of Mannheim partner $D$-curves include those of Mannheim partner curves in some special cases. Furthermore, the relations between the geodesic curvature, the normal curvature and the geodesic torsion of these curves are given.

## REFERENCES

[1] Blaschke, W., *"Differential Geometrie and Geometrischke Grundlagen ven Einsteins Relativitasttheorie Dover"*, New York, (1945).
[2] Burke J. F., *"Bertrand Curves Associated with a Pair of Curves"*, Mathematics Magazine, Vol. 34, No. 1. (Sep. - Oct., 1960), pp. 60-62.
[3] Görgülü, E., Ozdamar, E., *"A generalizations of the Bertrand curves as general inclined curves in $E^n$ "*, Communications de la Fac. Sci. Uni. Ankara, Series A1, 35 (1986), 53-60.
[4] Hacisalihoğlu, H. H., *"Diferansiyel Geometri"*, İnönü Üniversitesi Fen-Edebiyat Fakültesi Yayınları No:2, (1983).
[5] Izumiya, S., Takeuchi, N., *"Special Curves and Ruled surfaces"*, Beiträge zur Algebra und Geometrie Contributions to Algebra and Geometry, Vo. 44, No. 1, 203-212, 2003.
[6] Izumiya, S., Takeuchi, N., *"Generic properties of helices and Bertrand curves"*, Journal of Geometry 74 (2002) 97–109.
[7] Kazaz, M., Onder, M., *"Mannheim offsets of timelike ruled surfaces in Minkowski 3-space"*, arXiv:0906.2077v3 [math.DG].
[8] Kazaz M., Uğurlu H. H., Onder M., *"Mannheim offsets of spacelike ruled surfaces in Minkowski 3-space"*, arXiv:0906.4660v2 [math.DG].
[9] Kocayiğit, H., *"Minkowski 3-Uzayında Time-like Asal Normalli Space-like Eğrilerin Frenet ve Darboux vektörleri"*, C.B.Ü. Fen Bilimleri Enstisüsü Yüksek Lisans Tezi, 2004.
[10] Kurnaz, M., *"Timelike Regle Yüzeylerin Bertrand Ofsetleri"*, C.B.Ü. Fen Bilimleri Enstitüsü Yüksek Lisans Tezi, 2004.
[11] Liu, H., Wang, F., *"Mannheim partner curves in 3-space"*, Journal of Geometry, vol. 88, no. 1-2, pp. 120-126, 2008.
[12] O'Neill, B., *"Elemantery Differential Geometry"* Academic Press Inc. New York, 1966.
[13] O'Neill, B., *"Semi-Riemannian Geometry with Applications to Relativity"*, Academic Press, London, (1983).
[14] Orbay, K., Kasap, E., Aydemir, I, *"Mannheim Offsets of Ruled Surfaces"*, Mathematical Problems in Engineering, Volume 2009, Article ID 160917.
[15] Orbay, K., Kasap, E., *"On Mannheim Partner Curves in $E^3$ "*, International Journal of Physical Sciences Vol. 4 (5), pp. 261-264, May, 2009
[16] Ravani, B., Ku, T. S., *"Bertrand Offsets of ruled and developable surfaces"*, Comp. Aided Geom. Design, (23), No. 2, (1991).
[17] Struik, D. J., Lectures on Classical Differential Geometry, 2[nd] ed. Addison Wesley, Dover, (1988).
[18] Uğurlu, H. H., Kocayiğit, H., *"The Frenet and Darboux Instantaneous Rotain Vectors of Curves on Time-Like Surface"*, Mathematical&Computational Applications, Vol. 1, No.2, pp.133-141 (1996).
[19] Uğurlu, H. H., Topal, A., *"Relation Between Daboux Instantaneous Rotain Vectors of Curves on a Time-Like Surfaces"*, Mathematical&Computational Applications, Vol. 1, No.2, pp.149-157 (1996).
[20] Uğurlu, H. H., *"On the Geometry of Timelike Surfaces"* Communication, Ankara University, Faculty of Sciences, Dept. of Math., Series Al, Vol. 46, pp. 211-223(1997).

[21] Wang, F., Liu, H., *"Mannheim partner curves in 3-Euclidean space",* Mathematics in Practice and Theory, vol. 37, no. 1, pp. 141-143, 2007.

[22] Whittemore, J. K., *"Bertrand curves and helices"*, Duke Math. J. Volume 6, Number 1 (1940), 235-245.